\documentclass{nyj}

\title{Common divisors of $a^n-1$ and $b^n-1$ over function fields}

\author{Joseph H. Silverman}  
\address{Mathematics Department, Box 1917, Brown University,
Providence, RI 02912 USA}
\email{jhs@math.brown.edu}

\keywords{greatest common divisor, function field}
\subjclass{11T55; 11R58; 11D61}

\theoremstyle{plain}
\newtheorem{theorem}{Theorem}

\theoremstyle{definition}

\theoremstyle{remark}
\newtheorem{remark}{Remark}
\newtheorem{example}{Example}

\newtheorem*{acknowledgement}{Acknowledgements}

\newenvironment{notation}[0]{%
  \begin{list}%
    {}%
    {\setlength{\itemindent}{0pt}
     \setlength{\labelwidth}{4\parindent}
     \setlength{\labelsep}{\parindent}
     \setlength{\leftmargin}{5\parindent}
     \setlength{\itemsep}{0pt}
     }%
   }%
  {\end{list}}

  {\end{list}}

\renewcommand{\a}{\alpha}

\newcommand{\e}{\epsilon}
\newcommand{\f}{\phi}

\newcommand{\Ecal}{{\mathcal E}}

\newcommand{\Scal}{{\mathcal S}}

\renewcommand{\AA}{\mathbb{A}}

\newcommand{\CC}{\mathbb{C}}
\newcommand{\FF}{\mathbb{F}}

\newcommand{\QQ}{\mathbb{Q}}

\newcommand{\ZZ}{\mathbb{Z}}

\newcommand{\LCM}{\operatorname{LCM}}
\newcommand{\LS}[2]{\genfrac(){}{}{#1}{#2}}  

\newcommand{\notdivide}{\nmid}
\renewcommand{\pmod}[1]{~(\operatorname{mod}~#1)}

\begin{document} 

\begin{abstract}  
Ailon and Rudnick have shown that if $a,b\in\CC[T]$ are
multiplicatively independent polynomials, then
\text{$\deg\bigl(\gcd(a^n-1,b^n-1)\bigr)$} is bounded for all~$n\ge1$. We
show that if instead~$a,b\in\FF[T]$ for a finite field~$\FF$ of
characteristic~$p$, then \text{$\deg\bigl(\gcd(a^n-1,b^n-1)\bigr)$} is larger
than~$Cn$ for a constant $C=C(a,b)>0$ and for infinitely many~$n$,
even if~$n$ is restricted in various reasonable ways (e.g.,
$p\notdivide n$).
\end{abstract} 

\maketitle
\tableofcontents


\section*{Introduction}

Let~$a$ and~$b$ be positive integers that are multiplicatively
independent in~$\QQ^*$ and let $\e>0$.  Bugeaud, Corvaja, and
Zannier~\cite{BCZ} recently showed that there is an $n_0=n_0(a,b,\e)$
so that
\begin{equation}
  \label{equation:BCZ}
  \gcd(a^n-1,b^n-1) \le 2^{\e n}
  \qquad\text{for all $n\ge n_0$.}
\end{equation}
In other words, $a^n-1$ and $b^n-1$ cannot share a common factor of
significant size. Although elementary to state, the proof requires deep
tools from Diophantine analysis, specifically Schmidt's subspace
theorem~\cite{Schmidt}.
\par 
Ailon and Rudnick~\cite{AR} consider the analogous problem in which~$a$
and~$b$ are taken to be polynomials in~$\CC[T]$. They prove the stronger
result
\begin{equation}
  \label{equation:AR}
  \deg\bigl( \gcd(a^n-1,b^n-1)\bigr) \le C(a,b)
  \qquad\text{for all $n\ge1$.}
\end{equation}
It is natural to consider the situation when~$a$ and~$b$ are 
polynomials in~$\FF_q[T]$, where~$\FF_q$ is a finite field of
characteristic~$p$. In this case, some restriction on~$n$ is
certainly needed, since trivially
\[
  \gcd(a^{mp^k}-1,b^{mp^k}-1) = \gcd(a^m-1,b^m-1)^{p^k}.
\]
In this paper we will show that for $a,b\in\FF_q[T]$, even much
stronger restrictions on the allowable values of~$n$ do not allow
one to prove an estimate analogous to~\eqref{equation:BCZ}, much
less one as strong as~\eqref{equation:AR}.

\begin{acknowledgement}
The author would like to thank Gary Walsh for rekindling his interest
in arithmetic properties of divisibility sequences and for drawing his
attention to the papers~\cite{AR} and~\cite{BCZ}, Felipe Voloch for a
helpful discussion of Diophantine approximation in characteristic~$p$,
and Andrew Granville and the referee for several suggestions that
greatly improved this article.
\end{acknowledgement}


\section{A preliminary example}
\label{section:preliminaryexample}
As noted in the introduction, 
Ailon and Rudnick~\cite{AR} prove that if $a(T),b(T)\in\CC[T]$ are
nonconstant polynomials that are multiplicatively independent in~$\CC(T)$,
then
\[
  \deg\bigl( \gcd(a(T)^n-1,b(T)^n-1)\bigr) \le C(a,b)
  \qquad\text{for all $n\ge1$.}
\]
Suppose instead that $a(T),b(T)\in\FF_q[T]$, where~$\FF_q$ is a finite 
field of characteristic~$p$. It is natural to ask if the Ailon-Rudnick
estimate holds, at least if we require that $p\notdivide n$. 
As the following example shows, the answer is no.

\begin{example}
\label{example:example1}
Let $a(T)=T$ and $b(T)=T+1$, let~$Q=p^k$ be some power of~$p$,
and take \text{$n=Q-1$}. Then
\begin{multline*}
  b(T)^n-1 = (T+1)^{Q-1} - 1
  = \frac{(T+1)^Q-(T+1)}{T+1}
  = \frac{(T^Q+1)-(T+1)}{T+1} \\
  = \frac{T^Q-T}{T+1}
  = \frac{T(T^n-1)}{T+1}
  = \frac{T(a(T)^n-1)}{T+1}.
\end{multline*}
Hence
\begin{align*}
  \gcd(a(T)^n-1,b(T)^n-1) 
  &= \frac{\gcd((T+1)(a(T)^n-1),(T+1)(b(T)^n-1))}{T+1} \\
  &= \frac{\gcd((T+1)(a(T)^n-1),T(a(T)^n-1))}{T+1} \\
  &= \frac{a(T)^n-1}{T+1}.
\end{align*}
Therefore
\[
  \deg\bigl(\gcd(a(T)^n-1,b(T)^n-1) \bigr) = n-1.
\]
\end{example}

\section{Basic results on rational function fields over finite fields}
Before stating and proving a generalization of the example described
in Section~\ref{section:preliminaryexample}, we briefly recall some
basic arithmetic facts about the rational function field~$\FF_q(T)$.
We start with some notation:

\vbox{%
\begin{notation}
\item[$S_q$]
${}=\{\pi\in\FF_q[T] : \text{$\pi$ is monic and irreducible}\}$.
\item[$S_{q,N}$]
${}=\{\pi\in S_q : \deg(\pi) = N\}$.
\item[$S_{q,N}(\a,\mu)$]
${}=\{\pi\in S_{q,N} : \pi\equiv\a\pmod{\mu} \}$.
\item[$\Phi_q(\mu)$]
The function field analog of Euler's function~\cite[Chapter~1]{RosenFF},
\[
  \Phi_q(\mu) = \#\left(\frac{\FF_q[T]}{\mu\FF_q[T]}\right)^*
  = q^{\deg(\mu)} \prod_{\substack{\pi|\mu\\\pi\in S_q\\}} 
   \left(1-\frac{1}{q^{\deg(\pi)}}\right).
\]
\end{notation}
}

We are now ready to state three important theorems in the arithmetic
theory of (rational) function fields.

\begin{theorem}[Prime Number Theorem for Function Fields]
\label{theorem:primenumbertheorem}
\[
  \#S_{q,N} = \frac{q^N}{N}     + O\left(\frac{q^{N/2}}{N}\right).
\]
\end{theorem}
\begin{proof}
See \cite[Theorem~2.2]{RosenFF} for a proof. To see why this is the
analogue of the classical prime number theorem, notice that there
are~$q^N$ monic polynomials of degree~$N$ in~$\FF_q[T]$, so the fact
that~$\#S_{q,N}$ is asymptotic to $q^N/\log_q(q^N)$ is analogous to the fact
that~$\pi(X)$ is asymptotic to~$X/\log(X)$.
\end{proof}

\begin{theorem}[Primes in Arithmetic Progressions]
\label{theorem:dirichletstheorem}
Let $\a,\mu\in\FF_q[T]$ be  polynomials with $\gcd(\a,\mu)=1$. Then
\begin{equation}
  \label{equation:DirichletsTheorem}
  \#S_{q,N}(\a,\mu) = \frac{1}{\Phi_q(\mu)}\cdot\frac{q^N}{N}
    + O\left(\frac{q^{N/2}}{N}\right),
\end{equation}
\end{theorem}
\begin{proof}
See \cite[Theorem~4.8]{RosenFF} for a proof. Of course,
Theorem~\ref{theorem:primenumbertheorem} is the special case of
Theorem~\ref{theorem:dirichletstheorem} with $\a=0$ and $\mu=1$.
\end{proof}

Finally, we will need the following special case of the
general $r$-power reciprocity law in~$\FF_q[T]$. 

\begin{theorem}[$r$-Power Reciprocity]
\label{theorem:reciprocity}
Let $\pi\in S_q$, let $\mu\in\FF_q[T]$, and let~$r$
be an odd integer dividing~\text{$q-1$}. Then
\begin{equation}
  \label{equation:reciprocity}
  \pi\equiv1\pmod{\mu}~\text{{\upshape(}i.e., $\pi\in S_q(1,\mu)${\upshape)}}
   \Longrightarrow
  \text{$\mu$ is an $r^{\text{th}}$ power modulo $\pi$.}
\end{equation}
\end{theorem}
\begin{proof}
Let $\LS{\a}{\mu}_r\in\FF_q^*$ denote the $r^{\text{th}}$~power
residue symbol (\cite[Chapter~3]{RosenFF}), and let
$\{\a_1,\a_2,\ldots,\a_t\}$ be coset representatives for the
elements~$\a\in\FF_q[T]/\mu\FF_q[T]$ with the property that
\[
  \LS{\a}{\mu}_r = 1.
\]
Then one version~\cite[Proposition~3.6]{RosenFF} of the $r$-power
reciprocity law in~$\FF_q[T]$ says that for any monic irreducible
polynomial $\pi\in\FF_q[T]$,
\[
  \text{$\pi\equiv\a_i\pmod{\mu}$ for some $1\le i\le t$}
  \Longleftrightarrow
  \text{$\mu$ is an $r^{\text{th}}$ power modulo $\pi$.}
\]
(Note that our assumption that~$r$ is odd ensures that either
$(q-1)/r$ is even or else~$\FF_q$ has characteristic~2.)  In
particular, the implication~\eqref{equation:reciprocity} is an
immediate consequence of the fact that $\LS{1}{\mu}_r=1$ for every
modulus~$\mu$.
\end{proof}

\section{The main theorem}
\label{section:maintheorem}
Example~\ref{example:example1} shows that for particular
polynomials~$a(T)$ and~$b(T)$, the polynomial
\text{$\gcd(a(T)^n-1,b(T)^n-1)$} can be large when
\text{$n\equiv-1\pmod{q}$}. We first generalize this example to
arbitrary polynomials~$a(T)$ and~$b(T)$. We then consider more general
exponent values and show that it is unlikely that there is any
infinite ``natural'' set of exponents~$\Ecal$ with the property that
\[
  \bigl\{ n\in \Ecal : 
  \deg\bigl(\gcd(a(T)^n-1,b(T)^n-1) \bigr) \ge \e n\bigr\}
\]
is finite for every $\e>0$.

\begin{theorem}
\label{theorem:maintheorem}
Let $\FF_q$ be a finite field and let $a(T),b(T)\in\FF_q[T]$ be
nonconstant monic polynomials. Fix any power~$q^k$ of~$q$ and any
congruence class $n_0\in\ZZ/q^k\ZZ$. Then there is a positive constant
$c=c(a,b,q^k)>0$ such that
\[
  \deg\bigl(\gcd(a(T)^n-1,b(T)^n-1) \bigr) \ge cn 
  \qquad\text{for infinitely many $n\equiv n_0\pmod{q^k}$.}
\]
\end{theorem}

\begin{proof}
To illustrate the main ideas, we start with the special case $k=1$ and
$n_0=-1$, so as in Example~\ref{example:example1}, we look at
exponents~$n$ satisfying $n\equiv-1\pmod{q}$. More precisely, we will
take 
\[
   n=q^N-1\qquad\text{for $N=1,2,3,\ldots$.}
\] 
For all $\pi\in S_{q,N}$, the group
$\bigl(\FF_q[T]/(\pi)\bigr)^*\cong\FF_{q^N}^*$ has order~$n$, so as
long as $\pi\notdivide ab$, it follows that
\[
  a^n \equiv b^n \equiv 1 \pmod{\pi}.
\]
Hence for all sufficiently large~$N$, e.g. $N\ge\deg(ab)$, we have
\[
  \gcd(a^n-1,b^n-1)~\text{is divisible by}~\smash{\prod_{\pi\in S_{q,N}}} \pi,
\]
and hence
\[
  \deg\bigl(\gcd(a(T)^n-1,b(T)^n-1) \bigr) \ge \sum_{\pi\in S_{q,N}}\deg(\pi)
  = (\#S_{q,N})\cdot N.
\]
The ``Prime Number Theorem for Polynomials''
(Theorem~\ref{theorem:primenumbertheorem}) says that
\[
  \#S_{q,N} = q^N/N + O(q^{N/2}/N),
\]
so we find that
\[
  \deg\bigl(\gcd(a(T)^n-1,b(T)^n-1) \bigr)
  \ge q^N + O(q^{N/2})
  = n + O(\sqrt{n}\,).
\]
This completes the proof of the theorem for $n\equiv-1\pmod{q}$.  
\par 
In order to obtain more general exponents, we take~$n$ to have the
form \text{$(q^N-1)/r$} for a suitable choice of~$N$ and~$r$.  Then
\text{$\gcd(a^n-1,b^n-1)$} is divisible by primes~$\pi$ for which
both~$a$ and~$b$ are $r^{\text{th}}$~powers modulo~$\pi$. In order to
exploit this weaker condition, we will use the function field versions
of the power reciprocity law and Dirichlet's theorem on primes in
arithmetic progression.
\par 
For now, we  assume that $\gcd(q,n_0)=1$, since this is the most
interesting case. At the end of the proof we  briefly indicate what to
do if~$n_0$ is divisible by the characteristic.
Let $r\ge1$ be the smallest odd integer satisfying
\[
  r\cdot n_0 \equiv -1 \pmod{q^k},
\]
and let
\[
  Q = q^{k\f(r)},
\]
where $\f(\,\cdot\,)$ is the usual Euler phi function. We note that
\[
  Q \equiv 1 \pmod{r}.
\]
(If we were aiming for better constants, we could take~$Q$ to be
any power~$q^m$ with $m\ge k$ and $q^m\equiv1\pmod{r}$.)
\par
For each power~$Q^N$ of~$Q$, we let
\[
  n = n(Q^N) = \frac{Q^N-1}{r},
\]
and we observe that since $q^k|Q$, we have
\[
  r\cdot n \equiv -1\pmod{q^k},\quad\text{and hence}\quad
  n\equiv n_0\pmod{q^k}.
\] 
It remains to show that
$\deg\bigl(\gcd(a(T)^n-1,b(T)^n-1) \bigr) \ge cn$ for all $n=n(Q^N)$.
\par
Let
\[
  \ell(T) = \LCM[a(T),b(T)]\in\FF_q[T]
\]
be the (monic) least common multiple of~$a(T)$ and~$b(T)$.  We want to
use Theorem~\ref{theorem:reciprocity} to find primes for
which~$a(T)$ and~$b(T)$ are $r^{\text{th}}$~powers, but in order to
apply Theorem~\ref{theorem:reciprocity}, we need to work with a
sufficiently large base field. More precisely, we work in~$\FF_Q[T]$,
since our choice of~$Q$ ensures that the condition \text{$r|Q-1$} in
Theorem~\ref{theorem:reciprocity} is satisfied.
\par
We consider polynomials $\pi\in S_{Q,N}(1,\ell)$, i.e., monic
irreducible polynomials of degree~$N$ in $\FF_Q[T]$ satisfying
$\pi\equiv1\pmod{\ell}$. Then
\begin{align*}
  \pi\in S_{Q,N}(1,\ell)
  &\Longrightarrow 
  \pi\equiv1\pmod{a}\quad\text{and}\quad\pi\equiv1\pmod{b}\\
  &\Longrightarrow 
  a\equiv A^r\pmod{\pi}\quad\text{and}\quad b\equiv B^r\pmod{\pi} 
       \hbox to3em{\hfil}\\
  &\omit\hfill for some $A,B\in\FF_Q[T]$, 
                   from Theorem~\ref{theorem:reciprocity}, \\
  &\Longrightarrow 
  a^n \equiv (A^r)^{(Q^N-1)/r} = A^{Q^N-1} \equiv 1 \pmod{\pi}
    \quad\text{since $\deg\pi=N$,} \\
  &\omit\hfil
  and similarly  $b^n\equiv 1\pmod{\pi}$.
\end{align*}
This proves that $\gcd(a(T)^n-1,b(T)^n-1)$ is divisible by every polynomial
in $S_{Q,N}(1,\ell)$, so we obtain the lower bound
\begin{align*}
  \deg\bigl(\gcd(a(T)^n-1,b(T)^n-1) \bigr)
  &\ge \sum_{\pi\in S_{Q,N}(1,\ell)} \deg(\pi) \\
  &= \#S_{Q,N}(1,\ell)\cdot N \\
  &= \frac{Q^N}{\Phi_Q(\ell)} + O(Q^{N/2})
    \qquad\text{from Theorem~\ref{theorem:dirichletstheorem},}
\end{align*}
where recall that
\[
  \Phi_Q(\ell) = \#\left(\frac{\FF_Q[T]}{\ell\FF_Q[T]}\right)^*
  = Q^{\deg(\ell)} 
   \prod_{\substack{\pi|\ell\\ \pi\in\Scal_Q\\}}
    \left(1-\frac{1}{Q^{\deg(\pi)}}\right).
\]
Finally, using the definition $n=(Q^N-1)/r$, we see that
\begin{align*}
  \deg\bigl(\gcd(a(T)^n-1,b(T)^n-1) \bigr) 
    \ge \frac{r}{\Phi_Q(\ell)}&\cdot n + O(\sqrt{n})\\
  &\text{for all $n=\dfrac{Q^N-1}{r}$,}
\end{align*}
where the big-$O$ constant is independent of~$N$.
This gives the desired result, with an explicit value for~$c$.
\par
We are left to deal with the case that~$n_0$ is divisible by the
characteristic~$p$ of~$\FF_q$. Write $n_0=p^i\cdot n_1$ with
$p\notdivide n_1$. From the result proven above, we know that
\begin{equation}
  \label{equation:maintheorem1}
  \deg\bigl(\gcd(a(T)^n-1,b(T)^n-1) \bigr) \ge cn
  \quad\text{for infinitely many $n\equiv n_1\pmod{q^k}$.}
\end{equation}
For these values of~$n$, it follows 
\begin{align*}
  \deg\bigl(\gcd(a(T)^{p^in}-1,b(T)^{p^in}-1) \bigr) 
  &=  \deg\bigl(\gcd((a(T)^n-1)^{p^i},(b(T)^n-1)^{p^i}) \bigr)  \\
  &=  p^i\deg\bigl(\gcd(a(T)^n-1,b(T)^n-1) \bigr) \\
  &\ge cp^i n
\end{align*}
and that
\[
  p^i n \equiv p^i n_1 = n_0 \pmod{q^k}.
\]
Thus the values~$p^i n$ with~$n$ satisfying~\eqref{equation:maintheorem1}
show that Theorem~\ref{theorem:maintheorem}
is true when~$p|n_0$.
\end{proof}

\begin{remark}
The proof of Theorem~\ref{theorem:maintheorem} 
shows that it is true for any
\[
  c < \frac{1}{\Phi_Q(\LCM[a(T),b(T)])},
\]
where~$Q$ is a certain explicit, but potentially quite large, power
of~$q$.  More precisely, we can take $Q=q^m$ for some $m<2kq^k$, but
this huge number is undoubtedly far from optimal.
\end{remark}

\begin{remark}
\label{remark:trichotomy}
We conclude this note by again displaying the striking ``trichotomy''
in the values of~$\gcd(a^n-1,b^n-1)$:
\begin{align*}
\ZZ:   &\qquad\log\bigl(\gcd(a^n-1,b^n-1)\bigr)\le \e n.\\
\FF_q[T]:&\qquad\deg\bigl(\gcd(a^n-1,b^n-1\bigr))\ge c n.\\
\CC[T]:&\qquad\deg\bigl(\gcd(a^n-1,b^n-1)\bigr) \le c.
\end{align*}
For $\FF_q[T]$ and $\CC[T]$, these estimates are best possible, aside from
the delicate question of the value of the constants. However, 
the situation for~$\ZZ$ is less clear. Bugeaud, Corvaja, and 
Zannier~\cite[Remark~2 after Theorem~1]{BCZ} note that there is an 
absolute constant $c=c(a,b)>0$ so that
\[
  \log\bigl(\gcd(a^n-1,b^n-1)\bigr)\ge\frac{cn}{\log\log n}
\]
holds for infinitely many~$n$. It is reasonable to guess that this
lower bound is also an upper bound, with an appropriate larger choice
of~$c$. However, it appears to be an extremely difficult problem to 
prove any upper bound of the form
\[
  \log\bigl(\gcd(a^n-1,b^n-1)\bigr)\le f(n)
\]
with $f(n)$ satisfying $f(n)/n\to0$.
\end{remark}

\end{document}